\documentstyle[11pt,righttag,verbatim]{amsart}
\newtheorem{thm}{Theorem}[section]
\newtheorem{lemma}{Lemma}[section]
\newtheorem{prop}{Proposition}[section]
\newtheorem{cor}{Corollary}[section]
\newtheorem{defn}{Definition}[section]
\theoremstyle{definition}
\newtheorem{example}{Example}[section]
\theoremstyle{remark}

\def \mrn {{\Bbb R}^n}
\def \mr {{\Bbb R}}
\def \Natural {{\Bbb N}}
\def \Real {{\Bbb R}}

\def \eps {\epsilon}

\newcommand{\Id}{\operatorname{Id}}

\def \ha {\frac{1}{2}}

\def \st {:}
\def \p {\partial}

\def \sub {\subset}

\def \rao#1 {\frac{\p}{\p #1} #1}

\newcommand{\sgn}{\operatorname{sgn}}

\def \la {\langle}
\def \ra {\rangle}
\def \s {S}
\def \expbit {e^{i\lambda(\s  z'.\mu - \sqrt{1+\s^2} |z|)}}
\def \ampbit { a \left( \frac{1}{\s|z|}, \s, \mu\right) }
\def \ra {\rangle}
\newcommand{\SDiff}{\operatorname{Diff_{sc}^{2}}}
\newcommand{\Diff}{\operatorname{Diff}}
\def \sctb {\mbox{}^{sc}T^{*}(X)}
\def \sctbb {\sctb_{|\partial X}}
\def \schd {\mbox{}^{sc}\Omega^{1/2}}
\def \curdot {{\cal C}}
\def \legc {X \times \p X,\tilde{G};\schd}
\def \mcal { }

\newcommand{\sca}{\operatorname{sc}}
\newcommand{\spa}{\operatorname{span}}

\numberwithin{equation}{section}

  \title[Higher order scattering]{Higher order
 scattering on asymptotically Euclidean Manifolds}
   
\author{T. J. Christiansen}
 
\author {M. S. Joshi}


\begin{document}
\maketitle

\begin{abstract}We develop a scattering theory for perturbations of powers
of the Laplacian on asymptotically Euclidean manifolds.  The (absolute)
scattering matrix is shown to be a Fourier integral operator associated
to the geodesic flow at time $\pi$ on the boundary.  Furthermore, it 
is shown that on $\Real^n$ the asymptotics of
certain short-range perturbations of $\Delta^k$ can be recovered from the
scattering matrix at a finite number of energies.
\end{abstract}

\section{introduction}
 In this paper, we develop a scattering theory for
powers of the Laplacian on a class of manifolds
which includes perturbations of Euclidean space and apply
this theory to obtain new inverse scattering results
on Euclidean space.
This theory is a natural extension of the
work of Melrose \cite{sslaes}, who developed a theory of
scattering for the Laplacian on asymptotically Euclidean manifolds.
We show that the higher order scattering matrix has very similar
properties in this case. In particular, we show that the
Melrose-Zworski calculus of Legendrian distributions, \cite{smagai},
can be applied 
to construct the Poisson operator for the scattering problem and thus
deduce that the scattering matrix is a Fourier integral operator
associated to geodesic flow at time $\pi.$  For higher
order operators, this result appears to
be new for the class of perturbations we consider even for $\Real^n$.
This theory is then also
applied to extend the inverse results of Joshi and S\'{a} Barreto
on recovering the asymptotics of perturbations,
\cite{coulomb}, \cite{smagailess}, \cite{recpoten}, \cite{metric}, 
\cite{magnetic},  to
this higher order case.  

An asymptotically Euclidean manifold is a smooth manifold with
boundary $(X,\p X)$ which is equipped with a scattering metric. A
scattering metric is a smooth Riemannian metric, $g,$ on the interior of $X$
which blows up in a prescribed way at the boundary: there exists a
product decomposition close to the boundary, $p \mapsto (x,y) \in
 [0,\eps) \times \p X,$ such that $g$ takes the form
\begin{equation}
g = \frac{dx^2}{x^4} + \frac{h(x,y,dy)}{x^2}, \label{scatmet}
\end{equation}
with $h$ smooth on the closed space and $h_{|x=0}$ a non-degenerate
metric on $\p X.$ This is slightly different from Melrose's definition
in \cite{sslaes} but was shown to be equivalent in \cite{metric}. The
co-tensor, $h_{|x=0},$ is independent of the decomposition chosen and
thus we have a natural metric on $\p X.$ We shall assume throughout
that a product decomposition close to the boundary has been chosen and
fixed. 

It is important to realize that $\mrn$ with the Euclidean
metric is a special case of such a manifold.  
To see this, put $x = |z|^{-1}$ and $\omega
= z |z|^{-1};$ we then have
\begin{equation}
dz^{2} = \frac{dx^{2}}{x^4} + \frac{d\omega^2}{x^2},
\end{equation}
which also shows that the induced metric on the sphere at infinity is
the Euclidean metric on the unit sphere. 

Melrose showed that given a smooth function, $f,$ on $\p X,$ and
$\lambda \in \mr \setminus
\{0\}$ that there is a unique function $u,$ smooth on the interior
of $X$, of
the form
\begin{equation}
e^{i\lambda/x} x^{\frac{n-1}{2}} f_{+} + e^{-i\lambda/x} x^{\frac{n-1}{2}} f_{-},
\end{equation}
with $f_{\pm}$ smooth functions on $(X,\p X)$ and $f_{-}$ restricted
to the boundary equal to $f$ such that $(\Delta - \lambda^2) u=0.$ The
scattering matrix is then defined to be the map, $S(\lambda),$  on
$C^{\infty}(\p X),$ defined by
\begin{equation}
S(\lambda): f \mapsto {f_{+}}_{|\p X}.
\end{equation}
It was shown in \cite{sslaes} that $S(\lambda)$ extends to a unitary
operator on $L^{2}(\p X)$ with the density induced by $h.$ 
Melrose and Zworski (\cite{smagai}) studied the micro-local 
structure of this operator and showed that it is a zeroth order,
classical Fourier integral operator associated to geodesic flow at
time $\pi.$ In the special case that $h(x,y,dy) = h(0,y,dy) + {\cal
O}(x^{\infty}),$ Christiansen (\cite{t1}) and Parnovski (\cite{pa})
showed that for $\lambda<0$, modulo
smoothing,
\begin{equation}
S(\lambda) = ie^{i\pi \sqrt{\Delta_{\p X} +
\frac{(n-2)^{2}}{4}}}. \label{tanyasthm} 
\end{equation}

Here we develop analogous results for perturbations of powers of the
Laplacian. For $k\in \Natural$, 
we define a short range perturbation of $\Delta^{k}$
 to be a symmetric differential
operator of order $2k-1$ which close to the boundary
can be written in local coordinates in the form $x^2
P(x,y,x^2 D_x, x D_y)$ where $P(x,y,\tau,\eta)$ is smooth in $(x,y)$ and
is a polynomial of order $2k-1$ in $(\tau,\eta).$  In particular, a
real-valued, 
smooth function vanishing to second order at $\p X$ defines a short
range perturbation. 
  In Section \ref{s:bst} only, we allow certain 
pseudodifferential perturbations as well.

We prove
\begin{thm}
Let $V$ be a short range perturbation of $\Delta^{k}$ and let $\lambda
\in \mr \setminus \{0\}$.  Then, given $f \in C^{\infty}(\p X)$, there exists a
smooth function, $u,$  on $X^{0}$ such that $(\Delta^{k} + V -
\lambda^{2k}) u=0$ and $u$ is of the form
\begin{equation}
e^{i\lambda/x} x^{\frac{n-1}{2}} f_{+} + e^{-i\lambda/x}
x^{\frac{n-1}{2}} f_{-}, 
\end{equation} with $f_{\pm}$ smooth functions on $X$ such that 
$f_{-}$ restricted to the boundary is equal to $f.$ The function $u$
is unique modulo smooth functions vanishing to infinite order at $\p X.$
\end{thm}
The non-uniqueness here corresponds to the possibility of embedded
discrete spectrum. 
The scattering matrix $S(\lambda)$ can then be defined precisely as
before as the indeterminacy will not affect the lead term at the boundary,
and we show that it has a unitary extension. 

We also prove
\begin{thm} \label{smagaithm}
Let $V$ be a short range perturbation of $\Delta^{k}$ and let $\lambda
\in \mr \setminus\{0\}$.  Then $S(\lambda)$ is a zeroth order classical Fourier
integral operator associated to geodesic flow at time $\pi$ on $\p X.$
\end{thm}

It is interesting to compare the scattering matrix for different
powers and the following is an immediate corollary to our construction.
\begin{cor} \label{compdiffords}
Let $0< k_1 \leq k_2 \in {\Bbb N},$ and suppose $V_{j}$ is a short range
perturbation of $\Delta^{k_j},$ of the form $x^{l} p_{j}(x,y,x^2 D_x,
xD_y),$ with $p_j(x,y,\xi,\eta)$ a polynomial in $(\xi,\eta).$ Let
$S_{j}(\lambda)$ be the scattering matrix associated to
$\Delta^{k_j}+V_j - \lambda^{2k_j}.$  Then $S_{1}(\lambda)
- S_{2}(\lambda)$ is a Fourier integral operator of order $1-l.$
\end{cor}

It is therefore immediate that \eqref{tanyasthm} holds also for the higher
order scattering matrix when $V=0.$  

In the special case of $\mrn$ we study the problem of recovering
asymptotics of a perturbation from scattering data. As in
\cite{magnetic}, we need an {\it aradiality} condition to recover the 
perturbation. 
In fact, in \cite{magnetic} it was observed that recovery is not
possible without it.  We shall say a perturbation is aradial modulo
Schwartz functions if it is asymptotically equal to a sum 
\begin{equation}
\sum \limits_{l=2}^{\infty} \sum \limits_{\alpha} f_{\alpha,-l} D_{z}^{\alpha},
\end{equation}
with each term $\sum \limits_{\alpha} f_{\alpha,-l} D_{z}^{\alpha}$ of
the form $|z|^{-l}$ times a composition of vector fields tangent to
the sphere. The aradiality conditions allows any zeroth order
perturbation. 
\begin{thm}\label{thm:inv}
Let $V_1, V_2$ be
short range perturbations of $\Delta^{k}$ on $(\mrn, dz^{2})$ such that
$V_1 - V_2$ is an aradial 
differential 
operator of order $l.$ Let $S_{j}(\lambda)$ be the scattering matrix
associated to  $\Delta^{k} + V_j - \lambda^{2k},$ and suppose that  for
$l+1$ values of $\lambda>0,$ $S_1(\lambda) - S_2 (\lambda)$
is smoothing.  Then the coefficients of $V_1$ and $V_2$ agree modulo
Schwartz functions. 
\end{thm}
Note that we assume neither that $V_1, V_2$ are of order $l$ nor that
they are 
aradial.

Our approach to this
higher order scattering problem
is highly influenced by that of Melrose, \cite{sslaes},
and that of Melrose-Zworski, \cite{smagai}. In particular, we use
the scattering calculus developed by Melrose and used for the case $k=1$
to study the general case and use techniques similar to those of
\cite{sslaes} to establish the existence of the scattering matrix. To
establish the micro-local structure of the scattering matrix we
proceed as in \cite{smagai} to construct the Poisson operator for the
scattering problem as a Legendrian distribution associated to a pair
of intersecting Legendrian submanifolds using the calculus developed
there. 

To prove the inverse result, we choose to follow the approach of
\cite{smagailess} rather than \cite{recpoten} in order to increase the
readability of the paper for non-experts. In particular, we establish
our results for the case of $\mrn$ without explicitly using the Melrose-Zworski
Legendrian calculus. The reader may regard these proofs as a warm-up
for the construction of the Poisson operator in the general case. 
As in \cite{coulomb}, \cite{smagailess}, \cite{recpoten}, \cite{metric}, 
\cite{magnetic}, the proof proceeds by establishing that the principal
symbol of the difference of the scattering matrices determines and is
determined by a weighted integral of the lead term of the difference
of the perturbations over geodesics of length $\pi.$ The injectivity of
this transformation is then deduced by using some elementary calculus
and some deep results of Bailey and Eastwood, \cite{baileast}, on the
integral geometry of tensor fields on projective space.  

The main reference for higher order scattering on $\mrn$ is Chapter 14
of \cite{tome2} where a scattering theory for perturbations of a much
more general class of constant coefficient operators on $\mrn$ is
developed - this extends ideas developed by Agmon and H\"{o}rmander in
\cite{aghor} and Agmon in \cite{agm}.  The complementary problem of
studying the recovery of compactly supported perturbations of higher
order operators on $\mrn$ was studied by Liu in \cite{liu}.

The results of Bailey and Eastwood (\cite{baileast})
which we use in
section \ref{s:ip} to prove the inverse results 
 are restricted
to projective spaces and spheres.  Since in addition, we believe the
principal interest in this problem is in Euclidean
space, we restrict ourselves to studying that
important special case as the analysis is much more
accessible.
 The problem of recovering the asymptotics of metrics
for the $k=1$ case has been studied in \cite{metric}
 and it is likely that a similar result could
be proven in the higher order case. 

The question of recovering the entire perturbation
from the scattering matrix at fixed energy is an
interesting one but is still open even in the case
with $k=1.$

We are grateful to Mike Eastwood for helpful conversations and to the
Royal Society for the travel grant which made those conversations
possible. We also thank the London Mathematical Society for supporting
this collaborative research through its small grants scheme.  The 
first author is 
grateful for partial support from a University of Missouri S.R.F.

\section{Basic scattering theory of $\Delta^k+V$}\label{s:bst}

In this section, we construct solutions of $\Delta^k
+V-\lambda^{2k}$ having specified behaviour at the boundary, leading to
the definition of the scattering matrix.

Throughout, $\Delta$ denotes the Laplacian associated to a scattering
metric on a compact manifold $(X, \p X)$,
$k$ is a positive integer,
and $\lambda \in \Real \setminus \{0\}$.
We assume a product decomposition with boundary defining function $x$
has been fixed. We call the smooth functions vanishing to infinite
order at the boundary the Schwartz functions on $X.$ We use the
product decomposition to extend functions on the boundary smoothly
into $X$ by making them constant in the normal direction and cutting
off.

The results of this section are very closely related to the results of 
\cite{sslaes}.  We have tried to state the results in a manner accessible
to those unfamiliar with that paper, but in order to avoid repetition we 
omit proofs which follow essentially as in \cite{sslaes}, only giving some
 indication of how they need to be modified to work in this setting.
We recall some of the definitions and 
results of \cite{sslaes} but refer the reader
to the original for full details.

In this section only, we will allow a somewhat larger class of perturbations
of $\Delta^k$.  In order to describe them, we recall some notation.

Let $\cal{V}_{\sca}(X)$ be the space of all smooth
vector fields of finite length 
with respect to a scattering metric.  Near a point on 
the boundary, $x^2\frac{\partial}
{\partial x}$ and $x\frac{\partial}{\partial y_i}$, $i=1,...,n-1,$ form
a basis for $\cal{V}_{\sca}(X)$,
 where $y_i$ are coordinates on $\partial X$.  The set
of scattering differential operators of order $m$ is 
$$\Diff^m_{\sca}(X)=\spa_{0\leq j\leq m}(\cal{V}_{\sca}(X))^j.$$
The short range perturbations which we are allowing in the remainder of
this paper are elements of $x^2\Diff^{2k-1}_{\sca}(X)$.  

Most of the results of this section hold for a wider class of perturbations:
elements of $\Psi_{\sca}^{2k-1,2}(X)$, part of the 
(small) calculus of scattering pseudodifferential operators.  We refer the 
reader to \cite[Sections 4 and 5]{sslaes} for the full definition and some
properties of $\Psi^{m,l}_{\sca}(X)$.  However, we remark that 
$x^j\Diff^m_{\sca}(X)\subset \Psi^{m,j}_{\sca}(X)$ and $(\Delta -z)^{-1}
\in \Psi^{-2,0}_{\sca}(X)$ when $z \not \in [0,\infty)$ 
(\cite[Theorem 1]{sslaes}).  Moreover, if we consider the manifold
$(\Real^n,dz)$ (which 
by an earlier discussion under radial compactification becomes a compact
manifold $(\Real^n)_{rc}$
with scattering metric), the operators corresponding to 
$\Psi_{\sca}^{m,l}((\Real^n)_{rc})$ have Schwartz kernels of the form
$$A(z,z^{\prime})=\frac{1}{(2\pi)^n}\int e^{i(z-z^{\prime})\cdot \zeta}
a_L(z,\zeta)d\zeta$$
with $a_L$ satisfying
$$|D^{\alpha}_zD^{\beta}_{\zeta} a_L(z,\zeta)|\leq C_{\alpha \beta}
(1+|z|)^{-l-|\alpha|}(1+|\zeta|)^{m-|\beta|}
$$ (\cite[(4.1),(4.2)]{sslaes}).

\begin{defn}
We shall say a differential operator $V$ is a short range
perturbation of $\Delta^{k}$ if it is symmetric and if $V\in x^2
\Diff^{2k-1}_{\sca}(X)$.
We call a pseudodifferential operator a pseudodifferential short 
range perturbation of $\Delta^{k}$ if it is symmetric and it
is an element of $\Psi^{2k-1,2}_{\sca}(X)$.
\end{defn}

Throughout this section $V$ is a pseudodifferential short range perturbation
of $\Delta^k$.

We prove
\begin{prop} \label{prop:formexp}
Let $V$ be a pseudodifferential
short range perturbation of $\Delta^{k}$ and $\lambda \in
\mr \setminus\{0\}.$ Given $f \in C^{\infty}(\p X),$ there exists $u 
\in C^{\infty}(X),$ unique modulo Schwartz functions, such that $u_{|\p
X}=f$ and 
$ (\Delta^{k}+ V -\lambda^{2k}) (e^{-i\lambda/x}x^{\frac{n-1}{2}} u) $
is Schwartz. 
\end{prop}
Note we make no assumptions on the sign of $\lambda.$

This will follow easily once we have proven
\begin{lemma}\label{l:left}
Let $V$ be a pseudodifferential 
short range perturbation of $\Delta^{k}$ and $\lambda \in
\mr \setminus \{0\}.$ If $f \in C^{\infty}(\p X),$  then
\begin{equation*}
(\Delta^{k}+ V -\lambda^{2k}) (e^{-i\lambda/x}x^{\frac{n-1}{2}+\alpha} f) =
e^{-i\lambda/x} (
 k \lambda^{2k-1} C_{\alpha}
x^{\frac{n+1}{2}+\alpha} f + x^{\frac{n+3}{2}+\alpha}g)
\end{equation*}
with $g \in C^{\infty}(X),$ where $C_0 =0$ and $C_{\alpha} \neq 0$ for
$\alpha \neq 0.$ 
\end{lemma}
\begin{pf}
For $k=1,$ this follows as in \cite{sslaes},
using \cite[Lemma 8]{smagai} for the mapping properties of 
pseudodifferential short range perturbations. For general $k,$ one
simply iterates. 
\end{pf}

Now to prove Proposition \ref{prop:formexp}, we simply choose the lead term
to be $f$ and then repeatedly iterate to compute all the terms of the
Taylor series.  The result follows from Borel's lemma.  The uniqueness
follows from the fact that $C_{\alpha} \neq 0$ for $\alpha \neq 0.$

We will repeatedly use the operator \begin{equation}
Q=Q(\lambda)=\sum_{j=0}^{k-1}\lambda^{2j}\Delta^{k-j-1}
\end{equation}
since $\Delta ^k-\lambda^{2k}=Q(\Delta -\lambda^2)=(\Delta -\lambda^2)Q$. 
Using the same techniques as Theorem 1 of \cite{sslaes}, one can show 
that for
$\lambda \in \Real \setminus \{0\}$,
$Q(\lambda)^{-1}\in \Psi^{-2k+2,0}_{\sca}(X)$.

A key idea of many of the proofs is that 
\begin{equation}\label{eq:q}
\Delta^k-\lambda^{2k}+V
=Q(\Delta -\lambda^2+V^{\prime}),
\end{equation} where $V^{\prime}=
Q^{-1}V\in \Psi_{\sca}^{1,2}(X)$.  The importance of this is that the symbol
of $V^{\prime}$ vanishes to second order at the boundary.  As noted in 
Remark 3 of \cite{sslaes}, many of the results there hold if $\Delta $ is 
replaced by $\Delta +W$, when $W\in x^2C^{\infty}(X)$, because
they depend on the properties of the principal symbol and the boundary
symbol, which
are unchanged by the addition of such a $W$.  What we are doing
here is allowing a somewhat more general perturbation, but with the same
kind of decay at the boundary.  Thus, the results of Melrose's Propositions
9-11 hold if $\Delta$ is replaced by $\Delta +W$, with 
$W\in \Psi^{1,2}_{\sca}(X)$.

The following proposition is closely related to Proposition 11 of \cite{sslaes}
and follows from the modifications of Propositions 9-11 indicated above 
and from (\ref{eq:q}).
\begin{prop} \label{prop:uniquenes}
If $\lambda \in \mr \setminus \{0\}$, $u\in L^2_{sc}(X)$,
$V$ is a pseudodifferential short range perturbation of 
$\Delta^k$ and 
$(\Delta^{k}+V-\lambda^{2k}) 
u =0$, then $u$ is Schwartz.  
\end{prop} 

If $k=1$ and $V$ is a short range (differential) perturbation, then
$\Delta +V$ has no positive eigenvalues.
However, for (differential) short range
perturbations of $\Delta^k$, $k>1$, the
situation can be different.
\begin{prop}
There are (differential) short-range perturbations $V$ of 
$\Delta^k$, $k>1$, even, such that $\Delta^k+V$ has positive eigenvalues.
\end{prop}
\begin{pf}
We give two ways of constructing such examples.  Let $j\in \Natural$.

Using the 
min-max principle, for any asymptotically Euclidean manifold $X$ one can
find
$\tilde{V}\in x^2C^{\infty}(X)$ 
so that $\Delta +\tilde{V}$ has a negative eigenvalue $-\sigma$.  Then
$(\Delta +\tilde{V})^{2j}$ has a positive eigenvalue $\sigma^{2j}$.  It is 
possible to choose $\tilde{V}$ to be compactly supported and such that
$\Delta+\tilde{V}$ has many negative eigenvalues.

To construct a Schwartz potential $V$ such that $\Delta^{2j}+V$ has
an embedded eigenvalue, choose $\tau>0$ and
let $x$ be a globally defined boundary-defining function, 
$x\in C^{\infty}(X)$.
  Then $$(\Delta^{2j}-\tau^{4j})
e^{-\tau/x}x^{(n-1)/2}=\cal{O}(x^{(n+3)/2}e^{-\tau/x}).$$  Just as in Proposition \ref{prop:formexp} we can use Lemma
\ref{l:left}  to successively solve away the errors, resulting in 
a $u$ which satisfies
$$(\Delta^{2j}-\tau^{4j})u=g=\cal{O}(x^{\infty}e^{-\tau/x})$$
with $u=e^{-\tau/x}x^{(n-1)/2}(1+\cal{O}(x))$.  
In fact, by doing the asymptotic
summation judiciously, we can ensure that $x^{-(n-1)/2}e^{\tau/x}u$ is 
nowhere vanishing.  Set $V=-u^{-1}g$.  Then $\Delta^{2j}+V$ has eigenvalue 
$\tau^{4j} $ with eigenfunction $u$. 
\end{pf}

We continue with the results needed to define the scattering matrix.
\begin{thm} \label{thm:absorption}
If $\lambda \in \mr \setminus \{0\}$,
$V\in \Psi^{2k-1,2}_{\sca}(X)$ is symmetric,
 and $f$ is a Schwartz function on $X$,
orthogonal to the $L^2$ null space of $\Delta^k+V-\lambda^{2k}$
if there is any, then
there exists $u$  with $x^{-(n-1)/2}e^{-i\lambda/x}u \in C^{\infty}(X)$
such that 
$$(\Delta^{k}+V-\lambda^{2k}) u = f.$$  
\end{thm}  This is an analogue of \cite[Propositions 12 and 14]{sslaes}.
\begin{pf}  As in Proposition 14 of \cite{sslaes}, 
$$u=\lim _{t \downarrow 0} (\Delta^k+V-(\lambda+it)^{2k})^{-1}f$$
where the limit exists in $x^{-1/2-\delta}H^{\infty}_{\sca}(X)$ for any
$\delta >0$.  To show that the limit exists, let $u_t=(\Delta^k+V
-(\lambda+it)^{2k})^{-1}f$.  Then, since $f$ is orthogonal to the null
space of $\Delta^k+V-\lambda^{2k}$, so is $u_t$ for $t>0$.  Having made 
this observation, the proof follows as in \cite[Proposition 14]{sslaes},
 giving us the same additional
microlocal regularity as well.

To finish the proof, one uses the analog of Proposition 12 of \cite{sslaes}.
That is, if $\tilde{u}\in C^{-\infty}(X)$, 
$WF^{*,-1/2}_{\sca}(\tilde{u})\cap R_+(\lambda)
= \emptyset$, and $(\Delta + Q^{-1}V-\lambda^2)\tilde{u}\in \dot{C}^{\infty}
(X),$ 
then $e^{-i\lambda/x}x^{-(n-1)/2}\tilde{u} \in C^{\infty}(X)$.  
The proof follows just
as the proof of \cite[Proposition 12]{sslaes}, first noting that 
$[\tilde{\Delta}_0,Q^{-1}V]\in \Psi^{2,1}_{\sca}(X)$ and thus
$$[\tilde{\Delta}_0,Q^{-1}V]:x^{s}H^{\infty}_{\sca}(X)\rightarrow x^{s+1}
H^{\infty}_{\sca}(X).$$

The second observation that
allows the proof to proceed as in \cite[Prop. 12]{sslaes} 
is that if $\Diff^l_c(X)u\subset
H^{\infty,m}_{\sca}(X)$, then $\Diff_c^{l-l^{\prime}}(X)A
\Diff^{l^{\prime}}_c(X)u\subset
H^{\infty,m+s}_{\sca}(X)$ when $l^{\prime}\leq l$ is a nonnegative integer
and $A\in \Psi_{\sca}^{*,s}(X)$.  This then allows the inductive step in 
the proof of \cite[Proposition 12]{sslaes} to proceed just as it does there.
\end{pf}

We shall also use the following boundary pairing result, the analogue of
\cite[Proposition 13]{sslaes}.
\begin{prop} \label{prop:bp}
Suppose $V\in \Psi_{\sca}^{2k-1,2}(X)$ is symmetric,
$\lambda \in \Real \setminus \{0\}$,
 $u_i$, $i=1,2$ satisfies $(\Delta ^k+V-\lambda^{2k})u_i
=f_i$, $f_i$ is Schwartz, and $u_i=x^{(n-1)/2}(e^{i\lambda/x}a^+_i
+e^{-i\lambda/x}a^-_i)$, with $a^{\pm}_i\in C^{\infty}(X).$  Let 
$b^{\pm}_i=(a_i^{\pm})_{|\partial X}$.  Then
\begin{equation*}
2ik\lambda^{2k-1}\int _{\partial X}(b_1^+\overline{b_2^+}-b_1^
-\overline{b_2^-})
dh = \int _X u_1 \overline{f_2}-f_1\overline{u_2} dg.
\end{equation*}
\end{prop}
\begin{pf} Using $Q=Q(\lambda)=\sum_{j=0}^{k-1}\lambda^{2j}\Delta^{k-j-1}$,
we have
\begin{align}\label{eq:bp}& 
\int _X u_1 \overline{f_2}-f_1\overline{u_2} dg \nonumber \\ 
&= \int_X u_1\overline{(\Delta -\lambda^2+VQ^{-1})Qu_2}
-Q(\Delta-\lambda^2+Q^{-1}V)u_1\overline{u_2} dg \nonumber \\
& =\int_X u_1\overline{(\Delta -\lambda^2+VQ^{-1})Qu_2}
-(\Delta-\lambda^2+Q^{-1}V)u_1\overline{Qu_2} dg 
\end{align}
as $Q$ is self-adjoint and 
$(\Delta-\lambda^2+Q^{-1}V)u_1\in \dot{C}^{\infty}(X)$.  Recall that $
V$ is
self-adjoint as well. 
Choose a function $\phi \in C^{\infty}(\Real)$ such that $\phi(t)=0$ if 
$t<1$ and $\phi(t)=2$ if $t>2$. 
Then (\ref{eq:bp}) is equal to 
\begin{multline*}
\lim_{\epsilon \downarrow 0} \int_X 
\phi(x/\epsilon)\left(u_1\overline{(\Delta -\lambda^2+VQ^{-1})Qu_2}
-(\Delta-\lambda^2+Q^{-1}V)u_1\overline{Qu_2} \right) dg \\
=\lim_{\epsilon \downarrow 0} \int_X[\Delta 
+ Q^{-1}V,\phi(x/\epsilon)] u_1 \overline{Qu_2}dg.
\end{multline*}
Since $V$  is short-range, $\lim _{\epsilon \downarrow 0}
\int_X [ Q^{-1}V,\phi(x/\epsilon)]u_1 \overline{Qu_2}dg =0.$
To compute the remainder, we just use the corresponding results
from \cite[Proposition 13]{sslaes} and the fact that at the boundary
of $X$, 
\begin{equation*}
Qu_2= k\lambda^{2k-2}x^{(n-1)/2}e^{i\lambda/x}b_2^+ 
+k\lambda^{2k-2}x^{(n-1)/2}e^{-i\lambda/x}b_2^-
+\cal{O}(x^{(n+1)/2}).
\end{equation*}
\end{pf}

We will need
\begin{thm} \label{thm:existence}
If $\lambda \in \mr \setminus\{0\}$,
$V$ is a pseudodifferential short range perturbation of $\Delta^k$,
and $f \in C^{\infty}(\p X)$, then there
exists $f_{\pm} \in C^{\infty}(X)$ such that $u = x^{\frac{n-1}{2}}
\left( e^{i\lambda/x} f_{+} + e^{-i\lambda/x} f_{-} \right)$
satisfies
$$(\Delta^{k}+V-\lambda^{2k}) u=0 $$
and the restriction of $f_{-}$ to $\p X$ is $f.$
The function $u$ is unique
if $\lambda^{2k}$ is not an eigenvalue of $\Delta^k+V$ and is unique
up to the addition of a Schwartz function if $\lambda^{2k}$ is
an eigenvalue.
\end{thm}
\begin{pf}
Most of the proof follows from Proposition \ref{prop:formexp} and Theorem
\ref{thm:absorption}.  
We use Theorem \ref{thm:absorption} to solve away the error obtained by
constructing the formal expansion as in Proposition \ref{prop:formexp}. 
 We note 
that
since the error we wish to solve away is of the type $(\Delta^k+V-\lambda^2)g$,
where $g\in x^{-1/2-\epsilon}L^2_{sc}(X)$, and
since eigenfunctions are Schwartz we may 
still integrate by parts to obtain that the error is orthogonal to the
eigenspace.
Finally, suppose that there are two such $u$ satisfying the conditions
of the theorem.  Then their difference $v$ satisfies
$(\Delta^k+V-\lambda^{2k})v=0$ and $v=x^{(n-1)/2}e^{i\lambda/x}g_+
+x^{(n+1)/2}e^{-i\lambda/x}g_-$ with $g_{\pm}\in C^{\infty}(X).$  Then
the boundary pairing result of Proposition \ref{prop:bp} gives us that
$(g_+)_{|\partial X}=0$, and thus $v$ is in $L^2_{sc}(X)$.  Consequently,
$v$ is an eigenfunction and is thus Schwartz.
\end{pf}

We can therefore define the scattering matrix as in the usual case:
\begin{defn} \label{def:scatmatrix}
The scattering matrix is a map on $C^{\infty}(\p X)$ taking a function
$f$ to $f_{+}$ restricted to $\p X$ where $f,f_{+}$ are as in
Theorem \ref{thm:existence}. 
\end{defn}
Using Proposition \ref{prop:bp} and the uniqueness it implies, it is 
easy to check that the scattering matrix can
be extended to a unitary operator on $L^2(\partial X)$ which satisfies 
$S(\lambda)^{-1}=S(-\lambda)$.

We remark that if $V\in \Psi^{2k-1,2}_{\sca}(X)$ is symmetric,
then $\Delta^k+V$, initially
considered as an operator on the
Schwartz functions, has a unique self-adjoint extension to a domain in
$L^2_{\sca}(X)$. 
We end this section by describing the continuous part of its
spectral measure.  Let $P$ be the
Poisson operator whose existence is given by Theorem \ref{thm:existence}:
$$P(\lambda):C^{\infty}(\partial X)\ni f \mapsto u\in x^{(n-1)/2}e^{i\lambda/x}
C^{\infty}(X)+x^{(n-1)/2}e^{-i\lambda/x}
C^{\infty}(X)$$
where $u$ is the function in Theorem \ref{thm:existence}.  Then the continuous
part of the spectral
measure of $\Delta^k+V$ is given by
$$dE_c(\lambda)=\frac{1}{2\pi}P(\lambda)P^*(\lambda)d\lambda$$
where $\lambda^{2k}$ is the spectral variable and 
$\lambda \in [0,\infty)$.  This can be seen by first 
noting that if $V=0$ this follows from \cite[Lemma 2.2]{t1}
(see also \cite[Lemma 5.2]{h-v}) and the 
general case then follows as in \cite[Appendix to XI.6]{r-s}.  We note
that $\Delta^k+V$ has
 no singular continuous spectrum, though, as previously noted,
there may be discrete spectrum.

\section{The Poisson Operator for Euclidean space} \label{poissrn}

As an introduction to the general case, we use the techniques of 
\cite{smagailess} to construct the Poisson operator in the Euclidean 
case, where the construction can be done more explicitly.
Following \cite{smagailess} we look for a Poisson operator in the form
$e^{i\lambda z.\omega}a(z,\omega)$ with $a$ a polyhomogeneous symbol in $z$ and
$\omega$ is a smooth parameter. 

Let $\tilde{a}(z,\omega)$ be any polyhomogeneous symbol in $z$, with 
$\omega$ a parameter.  Then
\begin{equation}
\Delta(e^{i\lambda z.\omega} \tilde{a}) = (\lambda^2 \tilde{a} - 2i\lambda
\omega.\frac{\p \tilde{a}}{\p z} + \Delta \tilde{a})e^{i\lambda z.\omega}, 
\end{equation} 
so we conclude that 
\begin{equation}
\Delta^{k}(e^{i\lambda z.\omega} \tilde{a}) = 
(\lambda^{2k} \tilde{a} - 2ki\lambda^{2k-1}
\omega.\frac{\p \tilde{a}}{\p x} + b)e^{i\lambda z.\omega},
\end{equation} 
with $b$ a symbol two orders lower than $\tilde{a}.$ 

Now $V$ is a short
range perturbation of 
$\Delta^{k}$, which for $\mrn$ means that 
\begin{equation}
V = \sum \limits_{|\alpha| \leq 2k-1} f_{\alpha} D_{x}^{\alpha},
\end{equation}
with $f_{\alpha} \in S^{-2}_{cl}(\mrn).$ 
We conclude that it maps $e^{i\lambda z.\omega}S^{m}$ to $e^{i\lambda
z.\omega} S^{m-2}$ and
hence that if $\tilde{a} \in S_{phg}^{m}$ we have, 
\begin{equation}
(\Delta^{k}+V - \lambda^{2k})(e^{i\lambda z.\omega} \tilde{a}) 
= (- 2ki\lambda^{2k-1}
\omega.\frac{\p \tilde{a}}{\p z} + c)e^{i\lambda z.\omega},
\end{equation} 
with $c \in S_{phg}^{m-2}.$

Returning to the Poisson operator, 
by taking the lead term of $a$ to be the constant function $1$ we have
an error in $S^{-2}_{phg}.$  To solve away this and subsequent
error terms, we proceed precisely as in \cite[Section 2]{smagailess}, using
a term in $S^{-j+1}_{phg}$
to solve away an error term in $S^{-j}_{phg}$ (modulo terms
in $S^{-j-1}_{phg}$).  We 
solve the transport
equations,
$$- 2ki\lambda^{2k-1}
\omega.\frac{\p b}{\p z} = d,$$ where $d$ is the error, 
along the geodesics on the unit sphere from $\omega$ to $-\omega.$ The
only difference from the $k=1$ case of \cite{smagailess} is a factor of 
$\lambda^{2(k-1)}.$ As before, the transport equation degenerates on
approach to $-\omega$ and in particular blows up like $(\pi -s)^{-r}$
as $s,$ the geodesic distance from $\omega,$ tends to $\pi,$ when
solving for the term of homogeneity of order $-r.$

So as before we need a different ansatz close to $-\omega.$ We can
locally in $\omega$ rotate our coordinate system so that $\omega$ is
the north pole. 
So close to the south pole, we look for the Poisson operator in the form
\begin{equation}\label{eq:a2}
\int \limits_{0}^{\infty} \int \fracwithdelims(){1}{\s|z|}^{\gamma} \s^{\alpha}
\expbit 
\ampbit d\s d\mu,
\end{equation}
with $a(t,\s,\mu)$
 a smooth function compactly supported on $[0,\eps) \times [0,\eps)
\times S^{n-2}.$ We denote the class of functions that can be written
in this form plus a Schwartz error by $I^{\gamma,\alpha}.$ We recall
from \cite{smagailess} that this class is asymptotically complete in
$\gamma.$  It follows from a stationary phase computation carried out
in \cite{smagailess} that away from the south pole this is 
 equivalent to the previous ansatz
with a symbol of order $-\gamma+\frac{n-1}{2}.$

We recall from \cite{smagailess},
\begin{prop} \label{thm:disasymp}
If $u(z,\omega) \in I^{\gamma,\alpha}$ 
and $f\in C^{\infty}(\partial X \times \partial X)$ then 
$$ e^{i\lambda|z|} \int u(|z| \theta,\omega) f(\theta, \omega) d\theta d\omega$$ is a
smooth symbolic 
function in $|z|$ of order $-1-\alpha$ and its lead coefficient is
$|z|^{-\alpha-1} \la K, f
\ra$ where $K$ is the pull-back of the Schwartz kernel of a
pseudo-differential operator of order $\alpha-\gamma-(n-2)$ by the map
$\theta \mapsto
-\theta.$ The principal symbol of $K$ determines and is determined by the lead
term of the symbol, $a(t,\s,\mu),$ of $u$ as $\s \to 0+.$
\end{prop}

We also have as a special case of Proposition 15 from \cite{smagai},
\begin{prop}
If $u \in I^{\gamma,\alpha}$, then 
$$(\Delta^{k} - \lambda^{2k})u \in I^{\gamma+1,\alpha+1},$$
and $$Vu \in I^{\gamma+2,\alpha+2}.$$
\end{prop}

It is also important to note that,
\begin{lemma} \label{lemma:autoimprov}
If $u \in I^{\gamma,\frac{n-3}{2}}$ and $(\Delta^{k} + V -
\lambda^{2k})u \in I^{\gamma+j,\frac{n-1}{2}}$, then $(\Delta^{k} + V -
\lambda^{2k})u \in I^{\gamma+j,\frac{n+1}{2}}.$
\end{lemma}
This follows from a slight extension of the argument in the proof of
Lemma 3.2 in \cite{smagailess}. 

We also recall from \cite{smagailess} that 
\begin{prop} \label{prop:asymcom}
If $u \in I^{\infty,\alpha} = \bigcap \limits_{\gamma} I^{\gamma,\alpha}$,
then $u = e^{-i\lambda|z|}f(z),$ with $f$ a classical symbol of order
$-\alpha-1.$ 
\end{prop}

To carry out our construction we first use the
original ansatz, obtaining a symbol which blows up on approach to the
south pole.  Near 
the south pole, we use the second ansatz, (\ref{eq:a2}),
 taking $\gamma = -\frac{n-1}{2}, \alpha= \frac{n-3}{2}.$
The argument is then identical to the one in \cite{smagailess}. 
At the end of the construction, we obtain an approximate Poisson operator,
$\tilde{P},$ such that
$(\Delta^{k} + V - \lambda^{2k})\tilde{P}(\lambda)$ is
Schwartz. This error will in fact be orthogonal to the eigenspace at
energy $\lambda$ (if there is one) as the eigenfunctions are
necessarily Schwartz and the orthogonality
 follows from self-adjointness and integration by
parts. We can therefore apply the resolvent (Theorem \ref{thm:absorption})
and gain a term of the
from $e^{i\lambda|z|}f$ with $f$ smooth, solving away the error.

\section{The inverse problem}\label{s:ip}

In this section we prove Theorem \ref{thm:inv}.
Suppose we have two short range perturbations, $V_1, V_2,$  of
$\Delta^{k}$ and we want to compare their scattering matrices. In
particular, suppose 
\begin{equation}\label{eq:forcing}
V_1 - V_2 = \sum \limits_{|\alpha| \leq 2k-1}
a_{\alpha}(z)D_{z}^{\alpha},
\end{equation}
with $a_{\alpha} \in S^{-1-r}_{phg}(\mrn).$ 
If we carry out our construction for each $V_j,$ then the first $r$
terms of the construction with the first ansatz will be the same but
the term of homogeneity $-r$ will be different. In particular, the
forcing terms in the transport equations will differ by
\begin{equation}\label{eq:ply}
 W_{-r}=e^{-i\lambda z.\omega}(V_1 - V_2) e^{i\lambda z.\omega} = \sum
\limits_{|\alpha| \leq 2k-1} 
\omega^{\alpha}a_{\alpha,-r-1} \lambda^{|\alpha|},
\end{equation} 
where $a_{\alpha,-r-1}$ is the lead term of $a_{\alpha}.$
We therefore conclude that the terms of homogeneity $-r$ will differ
at $(\gamma(s),\omega)$ by
\begin{equation}
\frac{i |z|^{-r}}{2k\lambda^{2k-1}(\sin s)^{r}} \int \limits_{0}^{s}
W_{-r}(\gamma(s'),\omega)(\sin s')^{r-1} ds',
\end{equation}
where $\gamma$ is a geodesic on the sphere running from $\omega$ to
$-\omega.$ 
As we have shown that the lead singularity of the first ansatz as $s
\to \pi-$ is essentially the principal symbol of the scattering matrix, we
conclude by the same argument that the difference of the scattering
matrices will be of order $-r$ and that the principal symbol of the
difference will determine and be determined by,
\begin{equation}\label{eq:it}
\frac{1}{\lambda^{2k-1}} \int \limits_{0}^{\pi}
W_{-r}(\gamma(s),\omega)(\sin s')^{r-1} ds,
\end{equation}
for all geodesics $\gamma$ with $\omega$ equal to $\gamma(0).$ 

Note that as $W$ depends polynomially on $\lambda$, so does the
principal symbol.  In particular, if we know the scattering matrix for
$2k$ different values of $\lambda>0$, then we can separate out the
parts coming from the 
differing orders of $|\alpha|$ in the forcing term. More generally, if
we have a perturbation of order $l \leq 2k-1$, then we determine the
asymptotics from $l+1$ values. 
So in the sequel, assuming the perturbation is of order $l\leq 2k-1$,
we only consider forcing terms of the form
\begin{equation}
\sum \limits_{|\alpha| = l}
a_{\alpha}(z)D_{z}^{\alpha},
\end{equation}
where $l \leq 2k-1.$  To solve the inverse problem, then, we must show
that $a_{\alpha}(z)$, $|\alpha|=l$, can be recovered from knowledge of
the transform (\ref{eq:it}), where we replace the sums in (\ref{eq:forcing})
and (\ref{eq:ply}) by sums only over the terms with $|\alpha|=l$.  That
is, we need to show that the transform (\ref{eq:it}), a
weighted integral along
geodesics of length $\pi$, is invertible.
The injectivity of the transform  was proven in \cite{recpoten} for the case
$l=0.$  For higher order perturbations ($l>0$) the question is 
more subtle and requires the notion of {\it aradiality}.
We say a perturbation is {\it aradial} if it has no radial
component, that is, if the perturbation is a span of vector fields
tangent to the sphere. Clearly, all zeroth order perturbations are
aradial. The injectivity of the transform (\ref{eq:it})
for aradial first order
perturbations was shown in \cite{magnetic} and for second order
aradial perturbations in \cite{metric}. The impossibility of
recovering first order perturbations that are not aradial was also shown in
\cite{magnetic}. We now look at the general case.  

First we re-express the forcing term more invariantly. We regard the
perturbation as an $l$-form, $\mu=\sum \limits_{|\alpha|=l}
a_{\alpha,-r-1}(z) dz^{\alpha}.$ 
Since $dz^{\alpha}$ is symmetric, $\mu$ is symmetric.
The
aradiality means that $\mu$ makes sense as a form on the sphere and is
determined by its values as a map from the tangent space of the sphere
to $\mr.$ 

To re-express the forcing term we rotate so that $\omega =
(0,\dots,0,1)$ and $\gamma(s) = (0,\dots,0,\sin(s), \cos(s)).$ We then
have $\frac{d\gamma}{ds}(s) = (0,\dots,0,\cos(s), -\sin(s)).$ The
transform is then just $\int \limits_{0}^{\pi}
a_{(0,\dots,0,l)}(\gamma(s)) (\sin s)^{r-1} ds.$  We show that this is
equal to  
$$ (-1)^l \int \limits_{0}^{\pi} (\sin s)^{l+r-1}\left \la \mu,
\frac{d\gamma}{ds}(s)\right \ra ds$$
where $\left \la  \mu,
\frac{d\gamma}{ds}(s)\right \ra$ is the pairing of $\mu$ with
$\frac{d\gamma}{ds}\otimes\cdot \cdot\cdot \otimes \frac{d\gamma}{ds}$ ($l$ 
copies).
To prove this, observe that as $\gamma$ lies in a plane this is really
a two-dimensional question, so without loss of generality its enough to
take $n=2.$ Now $\mu$ is a symmetric $l-$form so 
$$\la \mu, v \ra = \sum \limits_{j=0}^{l} 
a_{(j,l-j)} v^{(j,l-j)}.$$
Thus 
$$\la \mu, \frac{d\gamma}{ds}(s) \ra = \sum \limits_{j=0}^{l} 
a_{(j,l-j)}(\gamma(s)) (-\sin(s))^{l-j} (\cos(s))^{j}.$$ 
It follows from aradiality 
and the fact that the perturbation is of order $l$
that 
$$a_{(j,l-j)} = \binom{l}{j}\left(\frac{-\cos(s)}{\sin(s)}\right)^j
 a_{(0,l)},$$ 
and thus we have 
$$(-\sin(s))^{l} \la \mu, \frac{d\gamma}{ds}(s) \ra = a_{(0,l)}\sum \limits_{j=0}^{l} \binom{l}{j}
(\cos(s))^{2j} (\sin(s))^{2l-2j},$$ 
which is of course equal to $a_{(0,l)}.$ This establishes the
equality. 

We now want to show that if $\mu$ is a symmetric $l-$form on the
sphere and 
\begin{equation}
\int \limits_{0}^{\pi} \la \mu(\gamma(s)), \frac{d\gamma}{ds}(s) \ra
(\sin(s))^{l+r-1} ds =0,
\end{equation}
for all geodesics $\gamma$ on the sphere, then $\mu=0.$ Since we know
this is true for any geodesic $\gamma,$ we have 
$$I_{l+r-1,\gamma,\alpha} = \int \limits_{0}^{\pi} \la \mu(\gamma(s+\alpha)),
\frac{d\gamma}{ds}(s+\alpha) \ra
(\sin(s))^{l+r-1} ds =0$$
for any $\alpha.$ Differentiating with respect to
$\alpha$ (see \cite{recpoten}), we deduce that $I_{j,\gamma,\alpha}=0$
implies that $I_{j-2,\gamma,\alpha}=0$ and thus, differentiating
repeatedly, if $l+r-1$ is
even and $I_{l+r-1,\gamma,\alpha} = 0$ then $\mu$ is even. Also, if $l+r-1$ is
odd and $I_{l+r-1,\gamma,\alpha} = 0$, then $\mu$ is odd. 

Since we can always reduce by two, it is enough to consider the cases
where $r$ is $1$ or $2.$
Considering $r=1,$ we have for any geodesic
$\gamma$ that 
$$\int \limits_{0}^{\pi} (\sin(s))^{l} \la \mu(\gamma(s)),
\frac{d\gamma}{ds}(s) \ra ds =0.$$
If we take a geodesic starting at $z_n=0$, then we deduce that
$$\int \limits_{0}^{\pi} \la (z_{n}^{l} \mu)(\gamma(s)),
\frac{d\gamma}{ds}(s) \ra ds =0.$$ By rotational invariance we have
that
$$\int \limits_{0}^{\pi} \left\la (p \mu)(\gamma(s)),
\frac{d\gamma}{ds}(s) \right\ra ds =0,$$
for all homogeneous polynomials, $p,$ of order $l.$ As $p\mu$ is even,
we can regard it as a symmetric tensor on projective space which is in
the kernel of the generalized x-ray transform on symmetric
$l-$tensors. Fortunately, the kernel of this operator has been
identified by Bailey and Eastwood, \cite{baileast}. They showed that
the kernel is precisely the symmetrized covariant derivatives of
symmetric $(l-1)$-tensors. Therefore, to show that
$\mu=0$ we need to show that if
$p\mu$ is a symmetrized covariant derivative for all homogeneous
polynomials $p$ of order $l$, then $\mu=0.$ 

We show this by showing that $\mu$ has to vanish at the north pole.  By
rotational invariance it will then follow that $\mu$ vanishes
everywhere.  For the case $r=2,$ we
apply the same arguments to $z_n \mu$ and the result will follow from
the case $r=1.$ 

We first prove
\begin{lemma}
Let $\mu$ be a symmetric co-tensor of order $l$ on $\mr^{n-1}$ such
that $p\mu$ is a symmetrized covariant derivative for every homogeneous
polynomial of order $l$.
Then $\mu$ vanishes at the origin.
\end{lemma}
\begin{pf}
To prove this lemma, we introduce some new notation. If $\alpha=(
\alpha_1, \dots, \alpha_{r}) \in \{ 1,\dots,n-1\}^{r}$ then
\begin{equation}
\p_{\alpha} = \frac{\p}{\p x_{\alpha_1}} \dots \frac{\p}{\p x_{\alpha_r}}.
\end{equation}
Note that $\p_{\alpha} \neq \p^{\alpha}$ even when both sides make
sense. We also put 
\begin{equation}
dx_{\alpha} = dx_{\alpha_1} \dots dx_{\alpha_r}.
\end{equation}
Let $\tilde{\alpha}_{t} =
(\alpha_1,\dots,\alpha_{t-1},\alpha_{t+1}, \dots \alpha_{t}).$  

The symmetrized covariant derivative, $\nabla_{s} \eta,$ of a tensor
$\eta$ is obtained by taking the usual covariant derivative and then
averaging over the symmetric group to make it symmetric. On $\Real^{n-1}$,
if the symmetric $l-1$ tensor
$\eta = \sum \phi_{\gamma} dx_{\gamma},$ and $\nabla_{s} \eta
= \sum \psi_{\alpha} dx_{\alpha},$ then
\begin{equation}
\psi_{\alpha} = \frac{1}{l}\sum \limits_{j=1}^{l} \p_{\alpha_{j}}
\phi_{\tilde{\alpha}_j}.
\end{equation}

Our proof of the lemma follows from the observation that the symmetric
$l-1$ tensors must satisfy certain PDEs. If $\alpha,\beta \in \{
1,\dots,n-1\}^{r},$ we define an exchange to be a map swapping certain
places in $\alpha$ with ones in $\beta.$ We include the case where no
swaps take place. The order of the exchange is the number of swaps.
The sign of the exchange will
be $-1$ to the power of the order. If the exchange is $e,$
we denote the sign of $e$ by $\sgn(e)$ and the new value of $\alpha$
by $e(\alpha,\beta)^{(1)}$ and of $\beta$ by $e(\alpha,\beta)^{(2)}.$
For example, if $e$ exchanges $\alpha_1$ with $\beta_1$ then 
\begin{align*}
e(\alpha,\beta)^{(1)} &= (\beta_1,\alpha_{2},\dots,\alpha_{r}) \\
e(\alpha,\beta)^{(2)} &= (\alpha_1,\beta_{2},\dots,\beta_{r}) \\
\sgn(e) &= -1 .
\end{align*}

If $\alpha, \beta \in  \{ 1,\dots,n-1 \}^{l},$ and $\mu=\sum \psi_{\gamma}
dx_{\gamma}$ is a covariant symmetrized derivative then 
\begin{equation}
\sum \limits_{e} \sgn(e) \p_{e(\alpha,\beta)^{(2)}}
\psi_{e(\alpha,\beta)^{(1)}} =0,
\end{equation}
where the sum is taken over all exchanges. To prove this one
substitutes the expression for a covariant derivative and observes
that the terms with no swaps will cancel with some of the terms from
swaps of order one but that the swaps of order one will then have
remainder terms which are canceled by swaps of order two and so
on. 

Suppose we have that every homogeneous polynomial of order $l$
times $\mu$ is a symmetrized covariant derivative. We show that each
term $\psi_{\alpha}$ vanishes at the origin. We consider two
cases. The first case is that $\alpha$ does not contain all possible
values from $\{1,\dots,n-1\},$ (which will always happen when
$l<n-1$).  Let $r$ be the value not taken. Let $\beta=(r,\dots,r).$ We
have that
\begin{equation}
\sum \limits_{e} \sgn(e) \p_{e(\alpha,\beta)^{(2)}}
(x_{\beta}\psi_{e(\alpha,\beta)^{(1)}}) =0 \label{eq:hompolyderiv}
\end{equation}
where $x_{\beta}=x_{\beta_1}x_{\beta_2}\cdot \cdot \cdot x_{\beta_l}$.
Upon evaluating at $x=0$, all terms will vanish except
$\p_{\beta}(x_{\beta} \psi_{\alpha})$, which will equal $ \psi_{\alpha}$,
which must therefore be zero.  

In the second case, $\alpha$ takes all values from $\{1,\dots,n-1\}.$
We let $\beta_{j}$ equal $2$ if $\alpha_{j}$ equals $1$ and  $1$
otherwise. As before, \eqref{eq:hompolyderiv} is satisfied. Upon
evaluation at zero, the only terms that will not vanish are those for
which $e(\alpha,\beta) = (\alpha,\beta)$
(up to the order of $\alpha_j$ and the order of 
$\beta_j$) and such $e$ will have
positive sign as an even number of swaps will be involved. We
therefore have that a positive multiple of $\psi_{\alpha}$ is zero and
thus that $\psi_{\alpha}$ is zero. 
\end{pf}

This completes the proof for $\mr^{n-1}$ with the
Euclidean metric. We,
however, wish to obtain a similar result for the sphere. Projecting the
upper hemi-sphere onto the plane, we obtain coordinates on the
sphere and the metric agrees with the Euclidean one at the north
pole (which corresponds to the origin).  The symmetrized covariant
derivative will then agree with the Euclidean one up to terms
vanishing there. As we are applying the result to polynomials of order
$l$ these additional terms do not affect the argument and so we
conclude that the result will also  hold for the sphere.

So we have shown that if the principal symbol vanishes then the lead
term of an {\it
aradial} perturbation vanishes also and thus by induction the
asymptotics of aradial perturbations are recoverable from fixed energy
scattering data.

\section{Review of Legendrian Distributions} \label{revleg}

In this section we review and rephrase the material we need from
\cite{sslaes} and \cite{smagai}. Here $X$ is a compact
manifold with boundary $\p X$ and $g$ is a scattering metric on $X$
and we have chosen a product decomposition of the form \eqref{scatmet}.
 Our account
is necessarily brief and we refer the reader to \cite{sslaes} and
\cite{smagai} for more details. 

There is a natural bundle over $X$ called the scattering cotangent
bundle which is denoted $\sctb.$ This is the dual to the bundle of smooth
vector fields of bounded length with respect to some (and hence all)
scattering metrics on $X.$ The restriction of $\sctb$ to $\p X$ is
denoted $\sctbb$ and carries a natural contact structure. If $y$ are
local coordinates on $\p X$ and $\mu$ are the corresponding dual coordinates,
then $(y,\mu,\tau)$ form local coordinates on $\sctbb,$ where $\tau$ is the
coefficient of $\frac{dx}{x^2}.$ 

We recall from Section \ref{s:bst} that 
 a
differential operator $P(x,y,xD_y,x^2 D_x)$ will be in
$\Psi^{m,k}_{sc}(X,\schd)$   
if it is of order $m$ and the total symbol as an operator in $xD_y,
x^2 D_x $ vanishes to $k^{th} $ order at the boundary. 
It then has a well-defined symbol at the boundary
$$j(P) = x^{k} p_k + x^{k+1} p_{k+1}, \; p_k , p_{k+1} \in
C^{\infty}(\mr \times T^{*}\p X).$$

\begin{defn}
An intersecting pair with  conic points is a subset, $\widetilde{W},$ of
$\sctbb$ which is a union of the closure of a smooth Legendrian
submanifold, $W$, and $W^{\#},$ which is a finite union of global sections
of the form $W^{\#}(\lambda_j) = \{ (y,0,\lambda_j) \},$ and contains 
$\overline{W}\setminus W.$ We also require $\overline{W}$ to have an at
most conic singularity at $\mu=0$; that is, it is smooth if polar
coordinates are introduced along $\overline{W}\setminus W.$
\end{defn}
The process of introducing polar coordinates along $\overline{W}\setminus W$ 
can be
given an invariant meaning and is then called blow-up. We denote the
blown-up manifold by $\widehat{W}.$ 

The metric $g$ induces a metric $h$ on the boundary as  nearby it is of the form 
$$ \tau^{2} + h'(y,\mu) + xg'$$
as a function on $\mbox{}^{sc}T^{*}X,$  we obtain $h'$ from $h$ via the
isomorphism 
$$\mu. \frac{dy}{x} \longmapsto \mu.dy.$$  

\begin{example}
For each $y' \in \p X$ and $0 \neq \lambda \in \mr,$ let
$G_{y'}(\lambda)$ be equal to the set of
$(\tau,y,\mu),$ such that $\tau^{2} + |\mu|^2 =\lambda^2, \mu \neq 0, $
and putting $\mu = |\mu| \hat{\mu},$ such that
\begin{gather}
\begin{gathered}
         \tau = |\lambda| \cos (s) \\
          |\mu| = |\lambda| \sin (s) \\
(y,\hat{\mu}) = \exp(s H_{\ha h})(y',\hat{\mu}')
\end{gathered} \label{pacoparam}
\end{gather}
where $s \in (0,\pi),$ $(y',\hat{\mu}') \in T^{*} \p X,$ and 
$h(y',\hat{\mu}')=1.$ Then $G_{y'}(\lambda) \cup \{ (\lambda,y,0) \}$ is an
intersecting pair with conic points. We denote this pair $\widetilde{G}(\lambda).$ This is the
pair in which we are interested here. The set $G^{\sharp}(\lambda) = 
\{ (-\lambda,y,0,y',0) \}$ is also important in our construction. Note
that $G^{\sharp}(\lambda)$  is the initial or incoming surface and
that $G^{\sharp}(-\lambda)$ is the outgoing surface. Note that in the
coordinates defined by \eqref{pacoparam}, $G^{\sharp}(\lambda)$ is $s
= \pi$ and $G^{\sharp}(-\lambda)$ is $s=0,$  when $\lambda$ is
positive. 
\end{example}
Associated with these intersecting pairs at each conic point is a
unique homogeneous Lagrangian submanifold 
$\Lambda(\widetilde{W}, \lambda_{i})$ 
of $T^{*}(\p X).$ For the pair $\widetilde{G}(\lambda)$ we are interested in,
this is precisely the relation of being $\pi$ apart along a lifted geodesic.
(See Proposition 4 of \cite{smagai}.) For simplicity, we shall henceforth take
 $\lambda$ to be positive. The $\lambda$ negative case is similar, 
or could be deduced from the positive case. 

Melrose and Zworski associated to any such intersecting pair a 
class of smooth functions whose asymptotics on approach to the boundary are
determined by symbols on the Legendrians. A symbol bundle over the
smooth Legendrian $W(\lambda)$ in the pair $\widetilde{W}$ can be defined and
is denoted $\hat{E}^{m,p}.$ The sections of this bundle are of the
form
$$a S^{p-m} |dx|^{m-n/4}$$
where $a$ is a smooth section of $C^{\infty}(\hat{W}; \Omega_{b}^{\ha}
\otimes M_{\hat{H}})$, $S$ is a defining function of the boundary
of $W,$ $M$ is the Maslov bundle, and $\Omega_{b}^{\ha}$ is the
$b-$half density bundle.
For $G$ above, one could take $S= \sin s.$ Melrose and
Zworski remove this singularity at the endpoints by rescaling but for
us it will be easier not to do so. 

\begin{prop} \label{symbolmap}
If $\widetilde{W}(\lambda)$ is an intersecting pair with conic points then
there is a class of smooth half-densities on $X^{o}$, denoted
$I^{m,p}_{sc}(X, \widetilde{W}),$  such that 
$\bigcap \limits_{m,p} I^{m,p}_{sc}(X, \widetilde{W})$ is equal to the
class of 
half-densities vanishing to infinite order at the boundary. There
exists a symbol map 
$$\hat{\sigma}_{sc,m,p}: I^{m,p}_{sc}(X, \widetilde{W}, \schd) \to
C^{\infty}(\hat{W}; \hat{E}^{m,p})$$
which gives a short exact sequence
$$ 0 \to I^{m+1,p}_{sc}(X,\widetilde{W}, \schd) \to
I^{m,p}_{sc}(X,\widetilde{W}, \schd)  \to C^{\infty}(\hat{W};
\hat{E}^{m,p}) \to 0.$$
\end{prop}
This is Proposition 12 from \cite{smagai}. The Legendrian
half-densities of order $m$ are
given locally away from the conic points by oscillatory integrals
\begin{equation}
 u = (2\pi)^{-\frac{n}{4} - \frac{k}{2}} \int e^{i\phi(y,u)}
a(x,y,u) x^{m-k/2+n/4}du,
\end{equation}
with $a$ smooth on $[0,\eps) \times U \times U'$ with $U,U'$ open and
$\phi$ parameterizes $W,$ that is
\begin{equation}
W= \{ (y,-\phi(y,u), d_{y} \phi) \st d_{u}\phi=0 \}.
\end{equation}
Near the conic singularity a more general form is required and we
refer the reader to \cite{smagai}. The order $m$ here is adjusted by
$-k/2+n/4$ and so is inconsistent with section \ref{poissrn}. We keep
this inconsistency as the second definition 
provides for better invariance properties but less clarity. 

An important related fact
we  need to know is how Legendrian distributions map
under scattering 
pseudo-differential operators. We recall Proposition 13 from
\cite{smagai}.
\begin{prop}\label{symbmapping}
Suppose $P \in \Psi^{l,k}_{sc}(X, \schd)$ has symbol $x^k p_k +
x^{k+1} p_{k+1}$ with respect to a product decomposition of $X$ near
$\p X,$ and suppose that $$W \sub \mbox{}^{sc}T_{\p X}^{*}(X)$$ is a
smooth Legendre submanifold. Then for any $m \in \mr,$
\begin{gather}
P: I^{m}_{sc}(X, W; \schd) \to I_{sc}^{m+k}(X,W; \schd) \\
\sigma_{sc,m+k}(Pu) = ({p_k}_{|G})\sigma_{sc,m}(u) \otimes |dx|^{k}.
\end{gather}
Furthermore, if $p_k$ vanishes identically on $W$ then 
$$P: I^{m}_{sc}(X,W; \schd) \to I^{m+k+1}_{sc}(X,W; \schd)$$
and
\begin{gather*}
\sigma_{sc,m+k+1} (Pu) = \\
\left( \frac{1}{i} \left( {\mcal L}_{V} +
\left( \ha(k+1) + m - \frac{n}{4} \right) \frac{\p p_k}{\p \tau} \right)+
{p_{k+1}} _{|W} \right) a \otimes |dx|^{m+k+1-\frac{n}{4} }
\end{gather*}
where $\sigma_{sc,m}(u) = a \otimes |dx|^{m-\frac{n}{4}}$ and $V$ is
the rescaled Hamiltonian vector field associated to $p_k.$ 
\end{prop}
We omit the definition of the rescaled Hamiltonian vector field but
recall that for $\Delta$ on the pair $G$ we are studying it is equal to
$$2 \lambda \sin s \frac{\p}{\p s}$$
in the semi-global coordinates given by \eqref{pacoparam}.

We also need two push-forward theorems, Propositions
16 and 17,  from
\cite{smagai}. They relate the singularities of the scattering matrix
to the asymptotics in small $x$ of the Poisson operator. Given a product 
decomposition near the boundary,
there is a natural pairing  
\begin{gather}
B :{\mcal C}^{-\infty}(X, \schd) \times C^{\infty}(\p X; \schd) \to
C^{-\infty}([0,\epsilon), \schd) \\
B(u,f) = x^{\frac{n-1}{2}} \int_{\p X} u(x,y) f(y). \label{pairing}
\end{gather}
\begin{prop} \label{pushforward}
For any intersecting pair of Legendre submanifolds with conic points,
$W,$ the partial pairing \eqref{pairing} gives a map
$$B: I^{m,p}_{sc}(X,\widetilde{W}; \schd) \times C^{\infty}(\p X; \schd)
\mapsto \sum \limits_j I^{p+\frac{n-1}{4}}( [0,\epsilon) ,
W'(\bar{\tau}_j; \schd))$$
where the $W'(\bar{\tau}_j) = \{ (0, -\tau_j dx/x^2 )\}$ are the
Legendre submanifolds corresponding to the components of $W^{\#}$ and 
$$B(u,f) = \sum \limits_{j} e^{- i \bar{\tau}/x} x^{p+ n/4}
Q^{0}_{\bar{\tau}_j}(u,f) \fracwithdelims||{dx}{x^2}^{\ha} +
O(x^{p+n/4+1})$$
with $$Q_{\bar{\tau}}^{0}(u) \in I^{p-m-\frac{n-1}{4}}_{phg}(\p X,
\Lambda(\widetilde{W}, \bar{\tau})),$$
and the principal symbol of $Q_{\bar{\tau}}^{0}(u)$ determines and is determined by the lead 
singularity of the principal symbol of $u$ on $W$ on approach to
$W'(\bar{\tau}_j).$ 
\end{prop}

When the Legendrian distribution is actually associated to a smooth Legendrian
submanifold the push-forward becomes much simpler and this simplifies the
construction of the Poisson operator. 

\begin{prop}
If $G$ is a smooth Legendre variety and $u \in I^{m}_{sc}(X,G', \schd)
$ near $\tau = \bar{\tau},$ then the distribution $Q^{0}_{\bar{\tau}}$
is a Dirac delta distribution.
\end{prop}

\section{The Poisson operator in the general case}

In this section, we apply the calculus reviewed in Section \ref{revleg} to
construct the Poisson operator and prove that the higher order
scattering matrix is indeed a zeroth order, classical, 
Fourier integral operator, proving 
Theorem \ref{smagaithm} and Corollary \ref{compdiffords}.  We shall refer
heavily to \cite{smagai}, Section 15 as our construction is a modification of
the one there. 

We assume a product decomposition of $X$ close to the boundary of the
form \eqref{scatmet}
 has been chosen
and is fixed throughout this section. We then have as in \cite{smagai} that 
$\Delta$, the intrinsic Laplacian acting on scattering half-densities on $X$, 
induces an operator 
$$\Delta_X \in \SDiff (X \times \p X, \schd(X\times \p X))$$
by 
$$ \Delta_{X} \left( u \fracwithdelims||{dx}{x^2}^{\ha}
\fracwithdelims||{dy}{x^{n-1}}^{\ha}
\fracwithdelims||{dy'}{x^{n-1}}^{\ha}\right) =
\Delta \left( u(.,y') \fracwithdelims||{dx}{x^2}^{\ha}
\fracwithdelims||{dy}{x^{n-1}}^{\ha}
\right) \fracwithdelims||{dy'}{x^{n-1}}^{\ha} ,$$
where $(x,y,y')$ is a point in $X \times \p X.$
Throughout this section $V$ will be a short range higher order
perturbation, that is a scattering differential operator of order $2k-1$
as a differential operator and order $2$ at the boundary. 
  
We recall from \cite{recpoten}, using the notation of (\ref{scatmet}) that
\begin{lemma} The symbol at the boundary of $\Delta_X$ is $p = p_0 +
xp_1$ with $p_0=\tau^2+h(0,y,\mu)$ and 
with $p_1$  equal to $ -i(n-1)\tau + c$ where $c=\frac{\partial}{\partial x}
h(x,y,\mu)_{|x=0}$ is quadratic in $\mu$.
\end{lemma}
Our choice of product decomposition ensures that 
there is no $\tau$ term in $c.$ We need to compute the symbols at the
boundary of $\Delta_{X}^{k} - \lambda^{2k} + V.$ The short range
assumption on $V$ ensures that it does not contribute - indeed this is
one reason 
why this definition of short range is appropriate. Now we can
decompose
$$ (\Delta^{k} - \lambda^{2k}) = Q(\lambda) (\Delta_{X} - \lambda^2),$$
with $$Q(\lambda) = \sum \limits_{j=0}^{k-1} \lambda^{2j}
\Delta_{X}^{k-1-j}.$$ We deduce that the lead symbol of $(\Delta_{X}^{k} -
\lambda^{2k})$ is $p_{0}^{k} - \lambda^{2k}$ and the second symbol
(subprincipal term) is equal to $k\lambda^{2(k-1)}$ times $p_1$ on the
zero set of the principal symbol. It also follows that the Hamiltonian
on the boundary of the lead term is just $k\lambda^{2(k-1)}$ times that
of $p_0$ on the zero set of the principal symbol. The fact that both
these terms have been changed simply by multiplication by
$k\lambda^{2(k-1)}$ ensures the simplicity of extending results from
the $k=1$ case to the higher order case. 

We also note the following lemma which is important in our construction to show
that the transport equations are solvable. 
\begin{lemma} \label{onebetter}
If $L \in I^{m,-\frac{1}{4}} (X \times \p X, \tilde{G}(\lambda),
\schd)$ is such that 
$$(\Delta_{X}^{k} - \lambda^{2k} + V) L \in I^{m+j,\frac{3}{4}}(X \times \p X,
\tilde{G}(\lambda),
\schd),$$ 
then 
$$(\Delta_{X}^{k} - \lambda^{2k} + V) L \in I^{m+j,\frac{7}{4}}(X \times \p X,
\tilde{G}(\lambda),
\schd).$$ 
\end{lemma}
This is a modification of Lemma 15 from \cite{smagai} and in fact, the
$k=1$ case, though not explicitly stated,
is essential to the construction there also.
\begin{pf}
The proof is no different from that of the special case in $\mrn$
(Lemma \ref{lemma:autoimprov}).
\end{pf}

\begin{prop} \label{poiscon}
For any $0 \neq \lambda \in \mr$ there exists $$K \in I^{m,p}(X\times \p
X, \tilde{G}(\lambda);\schd)$$ such that 
\begin{gather*}
(\Delta_{X}^{k} - \lambda^{2k} + V) K \in \curdot^{\infty}(X \times \p X;
\schd) \text{ and } \\ 
Q_{\lambda}^{0}(K) = \Id, \end{gather*} with
$$m = -\frac{2n-1}{4},\;  p = -\frac{1}{4},$$ 
and the principal symbol of $K$ on $G$ is 
$$C \sin(s)^{\frac{n-1}{2}} 
\frac{|ds|^{\ha}|dy|^{\ha} |d\hat{\mu}|^{\ha}}{(\sin s)^{\ha}}
 |dx|^{m- \frac{2n-1}{4}},$$
where $C(y,\hat{\mu})$ is a non-zero smooth function. 
\end{prop}
\begin{pf}
As in  \cite{smagai}, we first construct $K^{b} \in I^{m,p}(X\times \p
X, \tilde{G}(\lambda);\schd)$ such that 
\begin{gather}
(\Delta^{k}_{X} - \lambda^{2k} + V) K^{b} \in I_{sc}^{\frac{3}{4} }(X \times \p
X, G^{\sharp}(-\lambda)) \text{ and } \\
Q_{\lambda}^{0}(K^{b}) =\Id.
\end{gather}
We construct $K^{b}$ as an asymptotic sum of $$K_j \in
I_{sc}^{-\frac{2n-1}{4}+j, -\frac{1}{4}}(X \times \p X,
\tilde{G}(\lambda);
\schd).$$ 
We wish 
\begin{gather}
(\Delta^{k}_{X} - \lambda^{2k} + V) K_0 \in I_{sc}^{-\frac{2n-1}{4}  +2,
\frac{3}{4} } (X \times \p X, \tilde{G}(\lambda); \schd) \text{ and }
\\
\sigma_{0}(Q_{\lambda}^{0}(K_0)) = \sigma_{0}(\Id) 
\end{gather}
and then it is automatically in $$ I_{sc}^{-\frac{2n-1}{4} +2,
\frac{7}{4}} (X \times \p X, \tilde{G}(\lambda); \schd) $$ by Lemma
\ref{onebetter}. We also want 
$$
(\Delta_{X}^{k} - \lambda^{2k} +V) \left( \sum \limits_{l=0}^{j-1} K_l \right) \in   
I_{sc}^{-\frac{2n-1}{4}  +j+2,
\frac{3}{4}} (X \times \p X, \tilde{G}(\lambda); \schd) $$
and this, of course, implies that it will also be an element of $$I_{sc}^{-\frac{2n-1}{4}  +j+2,
\frac{7}{4} } (X \times \p X, \tilde{G}(\lambda); \schd).$$ 

Now near $G \cap G^{\sharp}(\lambda),$ where $G$ is smooth, we can as in
\cite{smagai} give an explicit construction and it is then  only necessary to
have that the principal symbol of $Q_{\lambda}^{0}(K_{0})$ is equal to $1$
 to ensure that 
$Q_{\lambda}^{0}(K^{b}) = \Id.$  We look for
$K_j$ of the form
$$x^{j} e^{i\lambda \phi(y,y')/x} a_{j}(x,y,y',\lambda) v, a_j \in
C^{\infty}(X \times \p X), $$ with $v$ a fixed scattering half-density and 
$\phi$ the cosine of the Riemannian distance from $y$ to $y' .$ Let $a_{j}^{'}$ be the restriction of $a_j$ to
$x=0.$  Taking geodesic normal coordinates, $y,$ about each $y'$  the
transport equations for $a_{j}^{'}$  is of the
form 
$$( y \cdot \p_{y} +j )a_{j}^{'} + b_j a^{'}_j = c_j \in C^{\infty}(X\times \p X)$$ 
near $y=0$ where $c_0 $ is identically zero and $b_j$ vanishes quadratically at
$y=0.$  So as in \cite{smagai}, the
terms $K_j$ exist sc-microlocally close to $G^{\sharp}(\lambda).$

We now need to continue each $K_j$ up to $G^{\sharp}(-\lambda).$  We do so by
solving transport equations for the principal symbols and iteratively solving
away the error. 

The principal symbol of $K_0,$ $ \sigma_{m}(K_0) ,$ is of the form
$$b \frac{|ds|^{\ha} |dy|^{\ha} |d\hat{\mu}|^{\ha}}{(\sin s)^{\ha}}
 |dx|^{m- \frac{2n-1}{4}}.$$ On the lifted geodesic $\beta(s)$ the
sub-principal term $c(\beta(s))= 2k\lambda^{2k-1} \sin(s) d(\beta(s))$ for
some smooth $d.$  
 From Proposition \ref{symbmapping}, the transport
equation for $b$ is 
\begin{equation*}
\frac{2k\lambda^{2k-1}}{i} \left( \sin(s) \frac{d}{ds} + \left( \frac{1-n}{2}\right) \cos(s) + i \sin(s) d(\beta(s))\right) b  =0
\end{equation*}
Writing $\tilde{b} = (e^{i \int d(\beta(s')) ds'} \sin(s)^{\frac{1-n}{2} }) b,$
 we thus have 
$$ \frac{d\tilde{b}}{ds}  =0.$$ 

This means that
\begin{equation}
b =C \sin(s)^{\frac{n-1}{2}} e^{-i \int d(\beta(s')) ds'}.
\end{equation}

As $s \to \pi-,$ that is near $G^{\sharp}(\lambda),$ this has a singularity of
the form $(\pi - s)^{\frac{n-1}{2}}$ and as $s \to 0+,$ of the form
$s^{\frac{n-1}{2}}.$  

As the order on $G^{\sharp}(\pm \lambda)$ is $ -\frac{1}{4},$ 
 the order $p
-m$  is equal  to the order of singularity and thus the solution
of the transport equation is a legitimate symbol and we can construct $K_0.$ 

Now by Lemma \ref{onebetter}, $$(\Delta^{k}_X + V - \lambda^{2k})(K_0) \in
I^{-\frac{2n-1}{4}+2,-\frac{1}{4}+2}(\legc),$$ and we look for $$K_1 \in
I^{-\frac{2n-1}{4}+1,-\frac{1}{4}}(\legc),$$ such that 
$$(\Delta^{k}_X + V - \lambda^{2k})(K_0+K_1) \in I^{-\frac{2n-1}{4}+3,-\frac{1}{4}+1}(\legc),$$
and thus by Lemma \ref{onebetter} is in
$I^{-\frac{2n-1}{4}+3,-\frac{1}{4}+2}(\legc).$ Letting 
the principal symbol of $K_1$ be $b_1 |dx|^{-\frac{2n-1}{4}+1-\frac{2n-1}{4}}$
 times the trivializing density above,
we obtain a transport equation; arguing as above, it becomes
\begin{multline*}
\sin(s)^{\frac{n-1}{2}}  e^{-i \int d(\beta(s')) ds'} \frac{d}{ds} e^{i \int d(\beta(s')) ds'}\left(
\sin(s)^{\frac{1-n}{2}+1} 
b_1 \right) \\= g(s) e^{-i \int d(\beta(s')) ds'}
\sin(s)^{\frac{n-1}{2}} ,\end{multline*} with $g(s)$
a smooth function on  $[0,\pi]$ (and smoothly depending on the suppressed
parameters).  Canceling, we obtain that 
$$ \frac{d}{ds} \left(e^{i \int d(\beta(s') ds'}
\sin(s)^{\frac{1-n}{2}+1} 
b_1 \right) = g(s),$$
which has a solution in the appropriate symbol class. The same argument, after
appropriately shifting indices, constructs all the terms $K_j.$ 

Asymptotically summing, we obtain $K^{b}$ such that 
$$(\Delta^{k}_X + V - \lambda^{2k}) K^{b} \in 
I^{7/4}(G^{\sharp}(-\lambda)).$$ These errors can now be removed by
an iterative construction of their Taylor series, cf Lemma 16 of \cite{smagai},
and we obtain $K$ as desired. 
\end{pf}

We have therefore constructed the Poisson operator modulo smooth terms
vanishing to infinite order at the boundary. We wish to remove this
error by applying the resolvent.  As before, this is possible even if there is 
embedded discrete spectrum,
as 
the error is orthogonal to the eigenspace at that
energy. This is seen by a simple integration by 
parts, since the elements of the eigenspace are Schwartz 
they are orthogonal to the image of smooth
functions of tempered growth. We have thus constructed $P(\lambda)$,
the Poisson operator for the problem, as a paired Legendrian
distribution. The remainder of the proof that $S(\lambda)$ is a
Fourier integral operator now follows as in \cite{smagai}, Proposition
19. 

To prove Corollary \ref{compdiffords}, one observes that if
$P_{1}(\lambda)$ is the Poisson operator for $\Delta - \lambda^2,$
then
\begin{equation}
(\Delta^{k} + V - \lambda^{2k}) P_1(\lambda) = Q(\lambda) (\Delta -
\lambda^2)P_1(\lambda) + VP_1(\lambda) = VP_1(\lambda).
\end{equation}
This means that the Poisson operator for $\Delta^{k} + V -
\lambda^{2k}$ can be constructed as a perturbation of $P_1(\lambda)$ and
as $VP_1(\lambda) \in I^{-\frac{2n-1}{4} + l, -\frac{1}{4} +l},$ the
scattering matrices will agree to order $1-l,$ and the corollary
follows.  Note, however, that even if one fixes the perturbation, the
scattering matrices will not agree to order more than $1-l$ as the
solutions of the transport equations will differ.

{\sc 
Department of Mathematics,
University of Missouri,
Columbia, Missouri 65211  U.S.A.}

{\em E-mail address:} {\tt tjc@@math.missouri.edu}

{\sc 
Darwin College, Cambridge CB3 9EU  U.K.}

{\sc Current Address:  NatWest Group Risk,
Level 9,
135 Bishopsgate,
London EC2M 3UR
U.K.}

{\em E-mail address:}
{\tt markjoshi@@alum.mit.edu}

\end{document}